\documentclass{article}
\newcommand{\R}{I\!\! R}

\newcommand{\C}{I\!\! C}

\newcommand{\Q}{I\!\! Q}

\begin{document}

Tsemo Aristide

Visitor, Department of mathematics

University of Toronto,

100 George Street
 tsemoaristide@hotmail.com

\medskip

\bigskip

\bigskip

\centerline{\bf GEOMETRIC STRUCTURES ON FIELDS.}

\bigskip

\bigskip

{\bf 0. Introduction.}

\bigskip

Let $X$ be a connected differentiable manifold,  $G$ a Lie group
which acts transitively on $X$, and which action verifies the
unique extension property (u.e.p): this means that two elements of
$G$ which coincide on an open set of $X$, coincide on $X$. A
geometric structure on $X$, or an $(X,G)-$structure on $X$, is a
differentiable manifold $M$ endowed with an atlas
$(\phi_i:U_i\rightarrow X)_{i\in I}$, where $\phi_i$ is a
diffeomorphism onto its image such that:
$$
\phi_i\circ{\phi_j}^{-1}_{\mid U_i\cap U_j}:\phi_j(U_i\cap
U_j)\longrightarrow \phi_i(U_i\cap U_j)
$$
is the restriction of an element $g_{ij}$ of $G$. The family
$g_{ij}$ satisfies the Chasles relation:
$$
g_{ij}g_{jk}=g_{ik}.
$$
The $(X,G)$ structure of $M$ pulls back to its universal cover
$\hat M$. This last structure is defined by a local
diffeomorphism:
$$
D_M:\hat M\longrightarrow X
$$
called the developing map, which gives rise to a representation
$$
h_M:\pi_1(M)\longrightarrow G
$$called the holonomy representation of the $(X,G)$ structure of
$M$.

The developing map can also be constructed as follows: we consider
the sheaf of local $(X,G)$ transformations from $M$ to $X$, we
denote by $F$ its Etale space, we have a map:
$$
F\longrightarrow X
$$
$$
[f]_x\longrightarrow f(x),
$$
where $x$ is an element of $M$ and $[f]_x$ an element of the fiber
of $x$, this application is well defined since $G$ satisfied the
unique extension property. Its restriction to a connected
component of $F$ is the developing map up to a  cover.

\bigskip

Examples of such structures are:

 $n-$dimensional affine manifolds
where $X$ is ${\R}^n$ and $G$ is $Aff({\R}^n)$, the group of
affine transformations of ${\R}^n$,

$n-$dimensional projective manifolds where $X$ is  the projective
space $P{\R}^n$, and $G$ is the group of projective
transformations $PGl(n,{\R})$, ect...

 For more information about those structures, see

 We remark that to develop a theory of $(X,G)$ manifolds, we need
 essentially $1-$homotopy theory. The $1-$homotopy of toposes is
 well understood. The goal of this paper is to generalize this
 theory to toposes to give applications to algebraic
 geometric and fields theory.

 \bigskip

{\bf 1. The fundamental group of a topos.}

\bigskip

The theory of the fundamental group of a topos is well understood,
it is proposed as an exercice in 4. p. 321 exercise 2.7.5.

\bigskip

{\bf Definition 1.1.}

Let $C$ be a category, a sieve on $X$ is a subclass $R$ of the
class $Ob(C)$ of objects of $C$ such that for every map
$m:Y'\rightarrow Y$ such that $Y$ is in $R$, $Y'$ is also in $R$.
Let $f:C'\rightarrow C$ be functor, we denote by $R^f$ the pulls
back of $R$ to $C'$, it is the sieve of $C'$ which elements are
objects which image by $f$ is in $R$.

We denote by $C_Y$ subcategory of objects over $Y$.

\bigskip

{\bf Defnition 1.2.}

A  topology on $C$ is an application which assigns to each object
$S$ a non empty subclass $J(S)$ of the class of sieves over $S$,
such that

For every map $f:T\rightarrow S$, and every sieve $R$ of $J(S)$,
$R^f$ is a sieve of $J(T)$, (here $f$ is considered as a morphism
between the categories $C_T$ and $C_S$.

For every object $S$ of $C$, every element $R$ of $J(S)$, and
every element $R'$ of $C_X$, $R'$ is an element of $J(S)$ if for
every object $f:T\rightarrow S$ of $R$, $R^f$ is an object of
$J(T)$.

The elements of $J(S)$ are called the raffinements of $S$.

\bigskip

{\bf Definition. 1.3.}

A presheaf of sets on $C$, is a contravariant functor from $C$ to
the category of sets.

A sheaf of sets on $C$, is a presheaf of sets on $C$ such that for
every raffinement $R$ of $S$, the map
$$
F(S){\longrightarrow}{F_{\mid R}}
$$
is bijective, where $F_{\mid R}$ is the presheaf defined on $R$ by
$F_{\mid R}(f)=F(T)$, where $f:T\rightarrow S$ is a map of $R$.

\bigskip

{\bf Example.}

We consider the category $e$ whose the class of object has one
element $x$, and such that the set of morphisms of $x$ is the
singleton $\{Id_x\}$. A sheaf on this category, is a contravariant
functor which sends $x$ to a set.

Let $C$ be a category whose set of objects is not empty, there
exists a projection functor $e_C:C\rightarrow e$, which assigns
$e$ to each object of $C$.

\bigskip

{\bf Definition 1.4.}

A constant sheaf, is a sheaf which factor through $e$: that is,
such that there exists a sheaf $F_e$ of $e$ such that $F=F_e\circ
e_C$.

We will say that the sheaf $F$ is locally constant, if and only if
there is a covering family $(X_i)_{i\in I}$, such that the
restriction of $F$ to the category over $X_i$ is a constant sheaf.

\bigskip

In the sequel, we will suppose that the category $C$ is a topos,
this means that $C$ is equivalent to the category of sheaves
defined on a standard $U-$sieve for a given universe $U$, or that
one of the following properties is satisfied:

(i) Endowed with its canonical topology, $C$ is a $U-$sieve on
which every sheaf is representable,

(ii) Endowed with its canonical topology, $C$ is a $U-$sieve such
that:

- The projective limits exist on $C$.

- Summand indexed by an element of $U$ exist, are disjoint, an
universal

- Equivalence relations are effective and universal on $C$.

\bigskip

{\bf Definition 1.5.}

A morphism $f:X\rightarrow Y$ of toposes is a functor
$f^{-1}:Y\rightarrow X$ such that: the presheaf
$f_*(F)(Y)=F(f^{-1}(Y))$ is a sheaf.

The left adjoint functor $f^*$ of $f_*$ commute with finite
projective limits.

\bigskip

{\bf Definition 1.6.}

Let $Z$ be the final object of a topos we will say that $C$ is
connected if we cannot find a covering family $\{S_1,S_2\}$ of $C$
such that $S_1\times_Z S_2$ represents the empty object.

\bigskip

{\bf Definitions 1.7.}

- A topological covering family $(X_i)_{i\in I}$ of $C$ is
connected if for every $X_i$, the topos $C_{X_i}$ is connected.

- A topos is locally connected if and only if for every covering
family $(X_i)_{i\in I}$, there is a connected  coverng family
$(Y_j)_{j\in J}$ such that for each $j$, there is a $k(j)$ such
that $Y_j$ is a subobject of $X_{k(j)}$.

\bigskip

{\bf Definition 1.8.}

Let $(X_i)_{i\in I}$ be a connected covering family of a topos, a
path of this family denoted by $(i_1,...,i_n)$, will be a family
of object $(X_1,...,X_n)$ such that $X_i\times_ZX_{i+1}$ is
different from the initial object, where $Z$ is the final object.

\bigskip

{\bf Proposition 1.9.}

A locally connected topos $C$ is connected if and only if there is
a topological covering family $(X_i)_{i\in I}$ such that for each
object $X$ of $C$, and for each $i$ of $I$, there is a path
$(i_1=i,...,i_n)$ such that $X_n\times_ZX$ is not the initial
object.

\medskip

{\bf Proof.}

Let $C$ be a locally connected and connected topos. Suppose that
there is an object $X_i$ of the family $(X_i)_{i\in I}$ and an
object $X$ of $C$, such that for every path $(i_i=i,...,i_n)$
$X_n\times_ZX$ is  the initial object. Let $C(X_i)$ be the set of
objects of $C$ such that for each element $Y$ of $C(X_i)$, there
is a path $(i_1=i,...,i_n)$ such that $X_n\times_ZY$ is not the
initial object. The complementary $U_i$ is not the empty, and the
final object is the direct summand of the set of object of $U_i$
and the set of objects of $C(X_i)$.

\medskip

We will call $C(X_i)$ the connected component of $X_i$. A locally
connected topos is the direct summand of its connected components.

In the sequel, topos considered will be supposed locally connected
and connected.

\medskip

We endow the set of paths associated to the topological covering
family $(X_i)_{i\in I}$ of the topos $C$ with the following
relation:

Two paths $x$ and $y$ will be said equivalent if and only if there
is a sequence of paths $z_1,...,z_n$, and a $k\leq n$ such that
$z_k=(i_1,..i_l,i_{l+1}=i_l,..,i_m)$ and
$z_{k+1}=(i_1,..,i_l,i_{l+2},..,i_m)$, with $x=z_1$ and $y=z_n$.

We denote by $Path((X_i)_{i\in I})$, the set of equivalence
classes of paths of $(X_i)_{i\in I}$.

In each equivalence class of a path $x$, we can find a
representant $\bar x=(j_1,...,j_k)$, such that $j_r\neq j_{r+1}$.

\medskip

Let $x=(i_1,...,i_n)$ and $y=(j_1,...,j_m)$ be two paths
representing the elements $\bar x$ and $\bar y$ of
$Path((X_i)_{i\in I})$ such that $i_n=j_1$. We associate to $\bar
x$ and $\bar y$ the element $\bar{x*y}$ which has
$(i_1,..,i_n,j_2,..,j_m)$ as a representant.

We can now define the groupoid $Gr((X_i)_{i\in I}$, whose objects
are the elements of $I$. A morphism between $i$ and $j$ is the
class of a path $(i_1,...,i_n)$ such that $i_1=i$ and $i_n=j$.

The inverse of the path $(i_1,...,i_n)$ is $(i_n,...,i_1)$.

The set of morphisms represented by the elements
$(i_1=i,...,i_n=i)$ is a group denoted $Aut(i)$.

\medskip

Consider now a locally constant sheaf $F$ defined on $C$,
$(X_i)_{i\in I}$ a locally constant connected topological covering
family, such that the restriction of $F$ to $X_i$ is a constant
sheaf. Such a family will be called a trivializing family.

\medskip

Let $X$, and $Y$, be two objects of $C$ such that the restriction
of $F$ to $X$ and $Y$ is constant. Recall $Z$ is the final object
of $C$. We suppose that $X\times_ZY$ is not the initial object.
The projections $p_x:X\times_ZY\rightarrow X$ and
$p_y:X\times_ZY\rightarrow Y$ gives rise to the isomorphisms of
sets: $g_x:F(X)\rightarrow F(X\times_ZY)$ and $g_y:F(Y)\rightarrow
F(X\times_ZY)$. We deduce the isomorphism
$g_x^{-1}g_y:F(Y)\rightarrow F(X)$.

We can apply this to the path $(i_1,...,i_n)$, we deduce for each
$i_k$ an isomorphism $g_{i_ki_{k+1}}:F(X_{i_{k+1}})\rightarrow
F(X_{i{k}})$, and a morphism
$g_{i_1i_n`}=g_{i_1i_2}..g_{i_{n-1}i_n}:F(X_{i_n})\rightarrow
F(X_{i_1})$.

It results a representation
$$
hol_F: Aut(i)\longrightarrow Aut(F(X_i))
$$
where $Aut(F(X_i))$ is the group of automorphisms of the sets
$F(X_i)$. We set $F_S(X_i)=hol_F(Aut(i))$.

\medskip

{\bf  Proposition 1.10.}

{\it Let $C$ be a locally connected and connected topos $C$, and
$F$ a sheaf on $C$, for each objects $X_i$ of $X_j$ of the
connected trivializing family $(X_i)_{i\in I}$, the groups
$hol_F(X_i)$ and $hol_F(X_j)$ are isomorphic. Moreover the class
of isomorphism of $hol_F(X_i)$ does not depend of the trivializing
family.}

\medskip

{\bf Proof.}

Let $X_i$ and $X_j$ be two elements of $(X_i)_{i\in I}$, there
exists a path $(i_1,...,i_n)$ between $X_i$ and $X_j$. This path
induces an isomorphism between $hol_F(X_i)$ and $hol_F(X_j)$.

Now, we show that the isomorphism class of $hol_F(X_i)$ does not
depend of the chosen trivializing family. Let $(Y_j)_{j\in J}$ be
another trivializing family. There exists another trivializing
family $(Z_k)_{k\in K}$, such that each $Z_k$ is a sub object of
an object $X_{i(k)}$ of $(X_i)_{i\in I}$ and $Y_{j(k)}$ of
$(Y_j)_{j\in J}$.

It suffices to show that $hol_F(X_i)$ and $hol_F(Z_k)$ are
isomorphic for each element $Z_k$ of $(Z_k)_{k\in K}$.

Let $(k_1=k,...,k_n=k)$ be a path of $(Z_k)_{k\in K}$, for each
$k_i$, we choose an element $X_{k_i}$ of $(X_i)_{i\in I}$ such
that $Z_{k_i}$ is a sub object of $X_{k_i}$, thus  we obtain a
morphism of group from $hol_F((Z_k)_{k\in K})$ to
$hol_F((X_i)_{i\in I})$.

Define now its inverse. Let $(i_1=i,...,i_n=i)$ be a path of
$(X_i)_{i\in I}$, for every $X_{i_l}$, there exists an object
$Z_{k_l}$ of the family $(Z_k)_{k\in K}$ which is a subobject of
such that $Z_{k_l}\times_ZX_{k_{l+1}}$ is different from the
initial object. There exists a path between a component of
$Z_{k_l}\times_ZX_{k_{l+1}}$ and $Z_{k_{l+1}}$ since $X_{k_{l+1}}$
is connected. We thus obtain a path of $(Z_k)_{k\in K}$ we define
the inverse of the previous morphism.

\medskip

We will denote by $(G_i)_{i\in I}$, the family of groups such that
for each $i$ there exists a sheaf $F_i$ on $C$, such that $G_i$ is
the holonomy of $F_i$.

Given two sheaves $F_i$ and $F_j$, a morphism of sheaves
$f_{ij}:F_i\rightarrow F_j$ induces a morphism between the group
$G_i$ and $G_j$. In this case, we will say that $G_i\leq G_j$ if
the morphism $f_{ij}$ is a surjective morphism .

We thus define a projective system of groups $(G_i)_{i\in I}$ with
projective completion is the projective fundamental group of the
topos $C$. We denote it by $pro\pi_1(C)$.

\medskip

For a locally connected topos, we can define the profundamental
group of each of its connected component.

\medskip

{\bf Definition 1.11.}

We will say that a locally connected, and connected topos $C$ is
simply connected if and only if each locally constant sheaf on $C$
is a constant sheaf, otherwise, this means that the profundamental
group of $C$ is trivial.

We will say that $C$ is a locally simply connected topos $C$ if
for each topological covering family $(X_i)_{i\in I}$ of $C$,
there exists a sub family $(Y_j)_{j\in J}$ such that $Y_j$ is
simply connected.

\medskip

{\bf Remark.}

\medskip

For a locally connected and locally simply connected topos $C$,
the profundamental group of $C$ is a group called the fundamental
group which has the following description:

Let $(X_i)_{i\in I}$ be a topological covering family of $C$, such
that for each $i$, $X_i$ is connected and simply connected.

We denote by $G$ the intersection of the kernel of all
representation $Aut(i)\rightarrow hol_F(X_i)$, the fundamental
group is the quotient of $Aut(i)$ by $G$.

\bigskip

{\bf 2. Geometric structures in topos theory.}

\bigskip

The main goal of this part is to extend the notion of geometric
structures to the notion of topos.

\medskip

{\bf Definition 2.1.}

Let $C$ be a topos, and $G$ a subgroup of automorphisms of $C$, we
will say that $G$ satisfies the unique extension property (uep) if
two elements of $G$ which coincide on the category over an object
of $C$, agree on $C$.

\medskip

In the sequel, topoi considered, will be locally connected topoi,
and we will assume that they have a finite number of connected
components.

\medskip

{\bf Definition 2.2.}

Let $C$ be a topos, and $G$ a group of automorphisms of $C$ whose
elements verify the unique extension property. We will say that a
topos $D$ is a $(C,G)$ topos if and only if there exists a
topological covering family $(X_i)_{i\in I}$ of $D$ such that for
each $X_i$, there exists a local isomorphism
$$
\phi_i: X_i\longrightarrow C
$$

this means that there exists an object $Y_i$ of $C$, such that
$\phi_i$ factor through an isomorphism $\phi_i': X_i\rightarrow
Y_i$.

Suppose that $X_i\times_ZX_j$ is different from the initial
object, (recall that $Z$ is the final object) then there exists
and element $g_{ij}$ of $G$ such that

$$
g_{ij}(\phi_{j})_{X_i\times_ZX_j}={\phi_i}_{X_i\times X_j}
$$
where $(\phi_i)_{X_i\times_ZX_j}$ (resp.
$(\phi_j)_{X_i\times_ZX_j}$) is the restriction of $\phi_i$ (resp.
$\phi_j$) to the category above $X_i\times_ZX_j$.

We have
$$
g_{jk}(\phi_k)_{X_i\times_ZX_j\times_ZX_k}=
(\phi_j)_{X_i\times_ZX_j\times_ZX_k},
$$

$$
g_{ij}(\phi_j)_{X_i\times_ZX_j\times_ZX_k}=
(\phi_i)_{X_i\times_ZX_j\times_ZX_k}.
$$
We deduce that
$$
g_{ij}g_{jk}(\phi_k)_{X_i\times_ZX_j\times_ZX_k}=
(\phi_i)_{X_i\times_ZX_j\times_ZX_k}.
$$
The unique extension property implies that
$$
g_{ij}g_{jk}=g_{ik}\leqno{ (1)}.
$$
The relation $(1)$ allows us to define a locally constant sheaf
$F$ on $D$, such that the restriction of $F$ to the category above
$X_i$ is the constant sheaf $G$, for each object $X$ of $D$,
$F(X)$ is the kernel of

$$
\prod F(X_i){\buildrel{\rightarrow}\over{\rightarrow}}\prod
F(X_i\times_XX_j).
$$

The transitions morphisms of $F$ are induced by the elements
$g_{ij}$ of $G$.

\medskip

Let $D_0$ be a connected component of $D$, we deduce from the
universal property of $pro\pi_1(D_0)$, a representation:
$$
hol_{D_0,C,G}:pro\pi_1(D_0)\longrightarrow G
$$
given by the projection
$$
pro\pi_1(D_0)\longrightarrow hol_F(X_i).
$$
This representation is called the holonomy representation of the
$(C,G)$ structure of $D$.

The representations deduced from $X_i$ and $X_j$ are conjugated,
so the conjugacy class of the holonomy representation does not
depend of $X_i$.

\medskip

{\bf Definition 2.3.}

Let $D$ and $D'$ be two $(C,G)$ topoi, we will say that a local
isomorphism $f:D\rightarrow D'$ is a $(C,G)$ morphism if and only
if
$$
\psi_{i(j)}f_{\mid X_i}=k_{ij}\phi_i,
$$
where $k_{ij}$ is an element of $G$, $((X_i)_{i\in I},\phi_i)$ and
$((Y_j)_{j\in J},\psi_j)$ are two topological covering families
which define respectively the $(C,G)$ structures of $D$ and $D'$,
and $f_{\mid X_i}$ is the restriction of $f$ to the category above
$X_i$. Moreover we suppose that $f_{\mid X_i}$ factors through an
sub object $Y_{i(j)}$ of $D'$. We thus endowed the class of
$(C,G)$ structures to a structure of category.

\bigskip

Since the existence of a universal cover of the topos $D$ is not
sure, we replace it by the sheaf $H_S$ of local $(C,G)$ morphisms
from $D$ to $C$, in order to define a developing map.

To $H_S$ we can associate the category whose objects are couples
$(X,s)$ where $X$ is an object of $D$ and $s$ is an element of
$H_S(X)$. A map between two objects $(X,s)$ and $(Y,s')$ is a
morphism $f:X\rightarrow Y$ such that $H_S(f)(s')=s$. This
category $H$ is a topos.

We can define the morphism:
$$
Dev:H\longrightarrow C
$$
defined on the category above $(X,s)$ by $s$.

The morphism is called the developing map of the $(C,G)$ structure
of $D$.

\bigskip

We consider now  $C$ (resp. $C'$) a topos $C$,  and a group $G$
(resp. $G'$) of automorphisms of $C$ (resp. $C'$) whose elements
verify the unique extension property. Moreover, we suppose that we
have a local isomorphism of topoi $\phi:C\rightarrow C'$, and a
morphism of groups $\Phi:G\rightarrow G'$ such that for each
element $g$ of $G$, we have $\phi\circ g=\Phi(g)\circ \phi$.

To each $(C,G)$ topos defined by the covering family $(X_i)_{i\in
I}$, we can associate the $(C',G')$ topos defined by the family
$((X_i)_{i\in I},\phi\circ\phi_i)$.

We thus define a functor $\phi_*$ from the category of $(C,G)$
topoi to the category of $(C',G')$ topoi.

\bigskip

{\bf 3. The space of $(C,G)$ structure.}

\medskip

We consider a topos $D$ endowed with a $(C,G)$ structure. Let
$((X_i)_{i\in I}$ be the covering family which defines this
structure. Let consider the topos $S((X_i)_{i\in I})$ whose final
object is obtained by gluing the product of topoi $X_i\times C$ by
the relation
$$
\tilde g_{ij} U_i\times_ZU_j\times C\longrightarrow
U_i\times_ZU_j\times C
$$
induced by the transitions functions $g_{ij}$. $S((X_i)_{i\in I}$
is called the structural topos of the $(C,G)$ structure of $D$.

The family of projections $p_i:U_i\times C\rightarrow D$ define a
 projection $p:S((X_i)_{i\in I})\rightarrow D$, since the transition functions
 induced the identity on the first factor.

 \medskip

 {\bf Definition 1.3.}

 A morphism $s:D\rightarrow S((X_i)_{i\in I})$ such that $p\circ
 s=Id$ will be called a transverse section of the topos
 $S((X_i)_{i\in I}$.

 A section $s:D\rightarrow S((X_i)_{i\in I})$ is transverse if and
 only if $p_i\circ s_i$ is the identity, where $s_i$ is the
 restriction of $s$ to the category over $X_i$.

 \medskip

 A $(C,G)$ topos on $D$ is defined by a transverse section $s_0$.
 Conversely, we have the fact:

 \medskip

 {\bf Proposition 3.2.}

 {\it A transverse section to the topos $S((X_i)_{i\in I})$
 defines a $(C,G)$ structure on $D$.}

 \medskip

 The space of $(C,G)$ structures on the topos $D$ can also be
 describe as the set of the following triples:

 A $(C,G)$ structure on the topos $D$, defined by the trivializing
 family $((X_i)_{i\in I}$,

 - we fixe a chart $(X_{i_0},\phi_{i_0})$

- an isomorphism of topoi $D\rightarrow D'$.

\medskip

We denote $S(D,C,G)$ the set of $(C,G)$ structures of $D$.

We have an application:
$$
S(D,C,G)\longrightarrow Hol(pro\pi_1(D),G)
$$
which associates to each element of $S(D,C,G)$ its holonomy
representation.

\bigskip

{\bf 4. Applications to algebraic geometry.}

\medskip

We consider a scheme $S$ define on a field $k$, $G$ a group of
automorphisms of $S$ which verifies the unique extension property.
We can suppose for instance that $S$ is irreducible, and $G$ a
subgroup of isomorphisms of $S$. It acts also on the Etale sieve
$Et_S$ of $S$. We can define a notion of $(S,G)$ schemes in the
category of schemes defined over $k$.

An $(S,G)$ etale scheme is a scheme $T$ such that its Etale sieve
$Et_T$ is endowed with and $(S,G)$ structure. This means that
there exists an etale covering $(U_i)_{i\in I}$ of $T$, etale
morphisms
$$
f_i: Et_{U_i}\longrightarrow Et_S
$$
  such that there
 exists a morphism $g_{ij}$ of $G$ such that:
 $$
 {f_i}_{\mid U_i\times_TU_j}=g_{ij}{f_j}_{\mid U_i\times_TU_j}.
 $$

Let $\hat T\rightarrow T$ be an etale covering.
 We
denote by $\pi_{\hat T}$ the group of Deck transformations of this
covering map. We can lift the $(S,G)$ structure on $\hat T$, by
setting $\hat f_i:U_i\times_T\hat T\rightarrow  S=f_i\circ
p_{U_i}$, where $p_{U_i}$ is the first projection $U_i\times_T\hat
T\rightarrow U_i$.

\medskip

 Suppose that the holonomy group of the $(S,G)$ structure of $T$
 is finite, its define a sheaf $F$ with finite fiber over $T$. We
 deduce from the theory of fundamental group  that
 there exists a finite cover $\hat T$ such that the pulls back of
 $F$ on $\hat T$ is trivial. We called $\hat T$ a holonomy cover
 of the $(S,G)$ structure of $T$.

 Let $(X_i)_{i\in I}$ be the covering family which defines the
 $(S,G)$ structure of $T$. The family $(X_i\times_T\hat T)_{i\in
 I}$ define also a $(S,G)$ structure on  $\hat T$. Since the pulls-back
 of $F$ on $\hat T$ is trivial, $(X_i\times_T\hat T)$ is in fact a
 Zariski open set of $\hat T$. Thus, there exists a map $Dev:\hat
 T\rightarrow \hat S$, where $\hat S$ is a finite Etale cover of
 $M$, such the restriction of $Dev$ to
 $X_i\times_T\hat T$ is the map $X_i\times_T\hat T$ which defines
 the $(S,G)$ structure of $\hat T$.

 \medskip

In the general case, the inductive limits of the family
$(U_i)_{i\in I}$ defines an Artin space $\hat T$, such that the
pull-back of $F$, on $\hat T$ is trivial. This allow  to define a
developing map, the developing map $\hat T\rightarrow S$.

\medskip

{\bf Definition 3.3.}

The $(S,G)$ structure on $T$ is complete if and only if the
developing map defined on the preceding Artin cover space  is a
surjective etale map.

 \bigskip

 Consider an $(S,G)$ structure defined on the scheme $T$, with
 finite holonomy group $h(T)$, we denote $\hat T$ the
  finite covering space corresponding to $h(T)$. For a fixed $\hat T$, the
 set of deformations of the $(S,G)$ structures is the set of
 fibered
 spaces $D(T,\hat T,S,G)$ which are quotient defined by the relations:
 $$
  \hat T\times S\longrightarrow \hat T\times S
 $$
 $$
 (x,y)\longrightarrow (x,hol(T)(\gamma)(y))
 $$
where $hol(T)$ is the  holonomy representation of an $(S,G)$
structure of $T$.

For such a bundle, we have an horizontal foliation ${\cal F}$
which is the pull forward of the foliation of $\hat T\times S$
whose leaves are the subschemes $\hat T\times y$ where is a point
of $S$.

A section $s$ transverse  to this foliation defines the $(S,G)$
structure of $T$, since it lifts to a section $\hat s:\hat
T\rightarrow \hat T\times S$.

If $T'$ is another finite covering space of $T$ such that there
exists an etale morphism $\hat T\rightarrow  T'$, we have a map
$D(T,T',S,G)\rightarrow D(T,\hat T,S,G)$. The inductive limit of
the space $D(T,\hat T,S,G)$ is the classifying  space of the
$(S,G)$ structure of $T$. We denote it by $D(T,S,G)$.

\medskip

Suppose now that the schemes are defined on ${\C}$, we endowed the
group $G$ with the compact-open topology induced by the analytic
topology of $S$. We can now study the deformation space of the
representation $\rho:hol(T)\rightarrow G$.

\medskip

{\bf Proposition 4.1.}

{\it Let $(\rho_t)_{t\in I}$ be a deformation of the
representation $\rho=\rho_0$ of the holonomy of an $(S,G)$
structure of $T$ where $I$ is an open interval, there exists an
open interval $J$ which contains $0$ such that for every $t$ in
$I$,  the groups $\rho_0(hol(T))$ and $\rho_t(hol(T))$ are
isomorphic.}

\medskip

{\bf Proof.}

Let $g_1$,...,$g_n$ be the elements of $hol(T)$. Consider the set
$X$ of elements of $S$ such that the cardinal of the family
$(g_1(x),...,g_n(x))$ is $n$. The set $X$ is a non empty set.
There exists open sets $U_1$,...,$U_n$ disjoint each other
containing respectively $g_1(x)$,...,$g_n(x)$ for the induced
analytic topology, an open interval $J$ contained in $I$, such
that $\rho_t(g_i(x))$ is an element of $V_i$ for $i=1,...,n$,
where $V_i$ is an open set contained in $U_i$, moreover if
$g_ig_k(x)$ is an element of $V_l$, then $\rho_t(g_ig_k(x))$ is
also an element of $V_l$ for $t$ in $J$. We deduce that the map
$$
\rho_t(hol(T))\rightarrow \rho(hol(T))
$$
$$
\rho_t(g_i) \longrightarrow g_i
$$
is an isomorphism.

\medskip

{\bf Theorem 4.2.}

{\it Let $S$ and $T$ be two ${\C}-$schemes, compact for the
analytic topology, then the image of the map $D(T,\hat
T,S,G)\rightarrow Hom(hol(T),G)$ is open.}

{\bf Proof.}

Under the conditions of the theorem, for each deformation
$(\rho_t)_{t\in I}$ of the holonomy representation of $T$, there
exists an interval $J$ contained in $I$ such that the groups
$\rho(hol(T))$ and $\rho_t(hol(T))$ are isomorphic. The flat
bundle $D(\hat T,T,S,G)$ and the flat bundle induced by the
representation $\rho_t$ for $t\in J$ are isomorphic. We identified
them. We denote by ${\cal F}_t$ the horizontal foliation of the
flat bundle induced by $\rho_t, t \in J$. For $t$ small, this
foliation remains transverse to the vertical foliation. It implies
that $\rho_t$ is the holonomy of a $(S,G)$ structure of $T$.

\medskip

{\bf The Beyli theorem.}

\medskip

Let ${\bf P}$ be the projective line defined over the rational
numbers, it has been shown by Beyli that finite covers of
$C_0={\bf P}-\{0,1,\infty\}$ are algebraic curves defined by
equations with algebraic coefficients. If we denote by $G$ the
group of automorphisms of $C_0$, we can restate this result as
follows:

\medskip

{\bf Theorem 4.3.}

{\it An algebraic curve has a complete $(C_0,G)$ structure which
holonomy is finite if and only if  it has a finite cover which is
an algebraic curve defined by algebraic coefficient.}

\medskip

The interest of the set of finite covers of $C_0$ is emphasizes by
the fact that the Galois group $Gal(\bar{\Q}/{\Q})$ acts naturally
on it by acting on coefficients of algebraic polynomials. Since
the structure of  $Gal(\bar{\Q}/{\Q})$ is not well-known, the
previous action is a tool which allow to study this group.

We will define a natural action of $Gal(\bar{\Q}/{\Q})$ on the set
of $(C_0,G)$ structures. This will give another tool to study this
group.

\medskip

Let $C_{\bar {\Q}}$ be the $\bar{\Q}-$projective line without
$3-$points, and $G_{\bar{\Q}}$ a subgroup of its automorphisms
group. Consider a curve $C$ defined over $\bar{\Q}$ endowed with
an $(C_{\bar{\Q}},G_{\bar{\Q}})$ structure. This structure is
defined by an etale covering family $(U_i\rightarrow C)$ of $C$,
and morphisms  $f_i:U_i\rightarrow V_i$, where $V_i\rightarrow
C_{\bar{\Q}}$ is an etale morphism. The group $Gal(\bar
{\Q}/{\Q})$ acts naturally on the space of curves defined over
$\bar{\Q}$, one of its element $\sigma$, sends a curve $C$ defined
by a polynomial $P$ which coefficients are algebraic, to the curve
$C^{\sigma}$ defined by $P^{\sigma}$. For an element $\sigma$ of
$Gal(\bar {\Q}/{\Q})$, we will denote by $C^{\sigma}$ the image of
$C$ by $\sigma$. The action of $\sigma$ induces a map
${f_i}^{\sigma}:{U_i}^{\sigma}\rightarrow {V_i}^{\sigma}$, we
define an $(C_{\bar{\Q}},G_{\bar{\Q}})$ structure on $C$. Remark
that we can have $C^{\sigma}=C$, but the
$(C_{\bar{\Q}},G_{\bar{\Q}})$ are not fixed by $\sigma$.

\medskip

{\bf The action of the Galois group $Gl(\bar{\Q}/{\Q})$ on the
fundamental groupoid of $C_0$.}

\medskip

We can define the topology $S_{C_0}$ generated  by the family of
finite covers maps $C\rightarrow C_0$. We have seen that a finite
cover $C$ of $C_0$ is an algebraic curve, define by polynomial
which coefficients are  algebraic. Let $\sigma$ be an element of
$Gl(\bar{\Q}/{\Q})$, $\sigma$ acts on finite covers of $C_0$.
Consider a path $(C_1,...,C_n)$ which represents an element of the
groupoid $Gr(C_0)$ where the index family is the set of finite
cover of $C_0$, we define the action of $\sigma$ on $Gr(C_0)$ by
$$
\sigma(C_1,..,C_n)=({C_1}^{\sigma},...,{C_n}^{\sigma})
$$
For each locally constant sheaf $F$ on $S_{C_0}$, we can define
the sheaf $F^{\sigma}$ defined by
$
 F^{\sigma}(C)=F(C^{\sigma})
 $

The action of $\sigma$ on $Gr(C_0)$ induces an automorphism
$\bar\sigma$ on the holonomy groupoid $G_F$ of $F$. This action is
compatible with morphisms between locally constant sheaves. We
thus deduce an action of $\sigma$ on the fundamental groupoid of
$S_{C_0}$. If $C$ is defined by polynomials with rational
coefficients, this action induces an action of
$Gal(\bar{\Q}/{\Q})$ on the fundamental group of $S_{C_0}$ in $C$.

\bigskip

{\bf The case of field theory.}

\medskip

We will now apply the previous study to field Theory. Let $k$ be a
field, $E$ an extension of $k$. We denote by $G$ the Galois group
of this extension. The unique extension property is satisfied by
the Galois group $G(E:k)$ acting on $Spec(E)$, endowed with its
Zariski topology, since the unique non empty set of $Spec(E)$, is
$Spec(E)$.

Recall that an etale cover of $Spec(E)$, is the spectrum of a
finite product of separable extensions $(k_i)_{i\in I}$ of $E$.
The Galois group of the etale cover is the product of the Galois
groups of the extensions $E\rightarrow k_i$.

 We want to classify finite products of finite separable extension
 of $k$ that we call $k-$separable algebras, which are endowed with $(E,G)$
structures. Let $F=k_1\times..\times k_n$ be such a $k-$separable
algebra, an etale covering of $Spec(F)$ is given by $Spec(F')$,
where $F'$ is a finite product of separable extensions of $k$
which is an extension of $F$. Let $k_S$, be the separable closure
of $k$. The family of separable extensions
$(H_i=(k^1_i\times..\times k^n_i)_{i\in I}$ of $F$, is a covering
family if and only if the $k_S$ is the inductive limit of
$(k^j_i)$. This is equivalent to saying that $Gal(k_S:F)$ is the
projective limit of the family $(Gal(k_S:k^j_i))$.

Thus the $(S,G)$ structure on $Spec(F)$ is defined by a family
$(F_i)_{i\in I}$ of finite product of separable extensions of $k$
such that $Gal(k_S:k^j_i)$ is the projective limit of the family
$Gal(F_S:k^j_i)$, and such that for each $i\in I$, there exists a
morphism
$$
D'_i: Et_{F_i}\rightarrow Et_{E},
$$
 a morphism $g_{ij}$ of $G$ which satisfies
$g'_{ij}{D'_j}_{Et_{Spec(F_i)\times_FSpec(F_j)}}=
{D'_i}_{Et_{Spec(F_i)\times_FSpec(F_j)}}$. The Chasles relation
$g'_{ij}g'_{jk}=g'_{ik}$ is satisfied.

\medskip

We will consider only morphisms induced by an etale map $D'_i:
Spec(F_i)\rightarrow Spec(E_i)$. Such a morphism is induced by an
separable extension

$$
D_i:E_i\longrightarrow F_i
$$

Thus the family $(g'_{ij})$ define a flat sheaf $S(E,G)$ on the
etale site of $Spec(F)$, with holonomy representation

$$
hol:Gal(F_S:F)\longrightarrow Gal(E:k)
$$

\medskip

{\bf Proposition.}

{\it Endow a field  $F$  with an $(S,G)$ structure defined by the
family $(F_i)_{i\in I}$, suppose that each $F-$algebra $F_i$ is a
field, then there exists a  separable extension  $F_{S,G}$ of $F$,
such that every field $F_i$ is a subfield of $F_{S,G}$, we denote
by $i$ the canonical map $F_i\rightarrow F_{S,G}$, there exists a
map $D:E\rightarrow F_{S,G}$, and a representation
$hol':Gal(F_{S,G}:F)\rightarrow Gal(E:k)$. We have $i\circ D=D_i$
The fields $F_{S,G}$ is universal in the following sense: If $K$,
is a field such that there exist a map $D_K:E\rightarrow K$, a map
$i_K:F_i\rightarrow K$ such that  $i_K\circ D_K=D_i$ then there
exist a map $f:F_{S,G}\rightarrow K$, such that $D_K=f\circ D$.}

\medskip

{\bf Proof.}

 The field $F_{S,G}$ is  the Artin cover which trivializes
the bundle defined by the holonomy of the $(E,G)$ structure of
$F$.

\medskip

{\bf Remark.}

An $(E,G)$ structures allow us to construct representations of
galois groups of separable extensions of $F$ to $k-$vector spaces.

\medskip

{\bf Definition.}

Let $F$  (resp.$F'$), be a  $k-$separable algebra endowed with the
respective $(E,G)$ structure defined by the etale covering
$(Spec(F_i))_{i\in I}$ (resp. $(Spec(F'_{i'})_{i'\in I'}$, and the
family of maps  $f_i:Spec(F_i) \rightarrow Spec(E_i)$, (resp.
$f'_i: Spec(F'_i)\rightarrow Spec(E'_i)$).

 An $(E,G)$ morphism between $F$ and $F'$ is an etale
 morphism $h:Spec(F)\rightarrow Spec(F')$ induced by
 a family of maps $h_i:Spec(F_i)\rightarrow Spec(F'_i)$, such that there exists a map
 $g_{ij}$ of $G$ with the property:
 $f'_i\circ h_i\circ u_i=g_{ij}\circ f'_j\circ h_j\circ u_j$, where $u_i$ is the
 canonical map $Spec(F_i)\times Spec(F_j)\rightarrow Spec(F_i)$,
  $h_i$ is the lift of $h$ to $Spec(F_i)$, and $g_{ij}$ is an element of
  $G$.

 We denote by $Aut(F,E,G)$ the
group of $(E,G)$ maps of $F$.

\medskip

{\bf The complete structures.}

\medskip

Recall that an $(E,G)$ structure on the field $F$ is said to be
complete if the developing map is a covering map.

\medskip

{\bf Proposition.}

{\it Endow a field $F$  with an $(E,G)$ structure, suppose that
the $(E,G)$ structure is defined by fields $(F_i)_{i\in I}$, the
structure is complete if and only if $F_{S,G}$ is a separable
extension of $E$.
 The developing map is an isomorphism if and only if
 and $E$ is isomorphic to the separable of $F$.}

\medskip

{\bf Proof.}

The $(E,G)$ structure of $F$ is defined by a family of maps
$E\rightarrow E_i\rightarrow F_i$ where $E_i$ and $F_i$ are
respectively separable extensions of $E$ and $F$. This structure
is complete if and only if $lim(F_i)=$the separable closure of $F$
is a separable extension of $E$. This implies that the developing
map is an isomorphism if and only if $F$ is separately closed and
$E$ is isomorphic to $F$.

\medskip

On the other hand we have the following:

\medskip

{\bf Proposition.}

{\it Suppose that  $F$ is a separable extension of the field $E$,
then $F$ is endowed with a complete $(E,G)$ structure.}

\medskip

Now suppose that the holonomy group $hol(Gal(F_s:F))$ is finite
and the structure is complete, we have seen that this implies that
$E$ is a separable extension of $F$. If the characteristic of $k$
is zero, and $G=Gal(E\mid k)$, every subfield $F$ of $E$, such
that $E$ is a finite extension of $F$ endows $E$ with  with a
complete $(F,G_F)$ structure because every extension of $k$ is
separable.

\medskip

We can restate the fundamental theorem of the Galois theory as
follows:

\medskip

{\bf Proposition.}

{\it Let $E$ be a Galois extension of $k$ with galois group $G$,
then the complete $(F,G_F)$ structures on $E$ which holonomies
groups are subgroup of $G$,  is one to one with the subgroups of
$G(E:k)$.}

\medskip

{\bf Operations on $(E,G)$ structures.}

\medskip

The change of  basis.

Let $E\rightarrow E'$ a morphism of $k-$fields, and $G$ and $G'$
be the Galois groups of $E$ and $E'$ over $k$. Suppose moreover
that $E'$ is separable over $E$, then to any $(E',G')$ defined on
the field $F$ by the family of $F-$algebras $(F'_i)_{i\in I}$, and
the family of maps $f:\phi_i:E'_i\rightarrow F'_i$, we can assign
an $(E,G)$ structure defined by the same data, since in this case
separable extensions of $E'$ are also separable extensions of $E$.

We thus obtain a functor
$$
f^*:S(E',G')\rightarrow S(E,G).
$$

Conversely, consider an $(E,G)$ structure defined on the field
$F$, by the family of $F-$separable algebras $(F_i)_{i\in I}$, and
the family of $E-$separable algebras $(E_i)_{i\in I}$, the maps
$\phi: E_i\rightarrow F_i$, gives rise to maps
$\phi'_i:E_iE'\rightarrow F_iE'$, remark that the fields $E_iE'$
and $F_iE'$ are respectively separable extensions of $E'$, and
$FE'$.

Let $g_{ij}: E_i\otimes E_j\rightarrow E_i\otimes E_j$ be the
coordinates change of the $(E,G)$ structure on $F$. The
coordinates change of the $(E',G')$ structure on $FE'$ are defined
by $h_{ij}:E_iE'\otimes E_jE'=(E_i\otimes E_j)E'\rightarrow
E_iE'\otimes E_jE'=(E_i\otimes E_j)E'$, where
$h_{ij}=g_{ij}Id_{E'}$.

 we obtain an $(E',G')$ structure on $FE'$, we have just define a
functor
$$
f_*:S(E,G)\longrightarrow S(E',G').
$$

\medskip

{\bf Proposition.}

{\it Let $s$ be an $(E',G')$ structure defined on the field $F$,
the $S(E,G)$ structure $f_*f^*(s)$ is isomorphic to $s$.}

\medskip

{\bf Proof.}

This result from the trivial fact that for each extension
$E'\rightarrow F'$, the field $E'F'$ is canonical isomorphic to
$F'$.

\medskip

{\bf Proposition.}

{\it The functor $f_*$ is the left adjoint of the functor $f^*$,
that is we have an isomorphism between $Hom_{(E',G')}(f_*(s),t)$
and $Hom_{(E,G)}(s,f^*(t))$ where $s$ and $t$ are respectively an
$(E,G)$ and an $(E',G')$ structure.}

\medskip

{\bf Proof.}

 The forgetful  functor $h'$ from the category $Ext(E')$ of extensions of
 $E'$  to the category $Ext(E)$ of extensions of $E$ is  right adjoint to
 the functor $h$ from $Ext(E)$ to $Ext(E')$ which assigns to $F$,
 the field $FE'$. This implies the proposition.

 \medskip

Let $F$ and $F'$ be two $k-$algebras,  endowed with the respective
$(E,G)$ structures $A$ and $B$ defined respectively by  the
families of separable extensions $(F_i)_{i\in I}$, and
$(F'_j)_{j\in J}$ of $F$ and $F'$. The $k-$algebras $E_i\otimes
E_j$ and $F_i\otimes F_j$ are respectively separable extensions of
$E$ and $F\otimes F'$. The maps $D_i:E_i\rightarrow F_i$, and
$D_j:E_j\rightarrow F_j$, define  a map $D_{i\otimes j}:E_i\otimes
E_j\rightarrow F_i\otimes F'_j$, thus it defines an $(E,G)$
structure denoted $A\otimes B$.

We have just endowed the category of $(E,G)$ $k-$algebras with a
tensor product.

We will denote by $I$ be the $(E,G)$ structure defined on $E$ by
the identity of the etale site of $E$.

\medskip

Suppose now that the field $k$ is ${\Q}$, and $E$ is a product of
finite extensions of ${\Q}$. Let $F$, be a ${\Q}-$algebra endowed
with an $(E,G)$ structure, defined by the respective covering
families $(F_i)_{i\in I}$, and $(E_i)_{i\in I}$ of $F$ and $E$.
Let $\sigma$ be an element of $Gal(\bar{\Q}/{\Q})$. We will denote
by $L^{\sigma}$ the image of a ${\Q}-$algebra  $L$ by $\sigma$.
The $(E,G)$ structure of $F$ is defined by a family of maps
$f_i:E_i\rightarrow F_i$, the action of $\sigma$ on $\bar {\Q}$,
induces an action ${f_i}^{\sigma}:{E_i}^{\sigma}\rightarrow
{F_i}^{\sigma}$. This action defined and $(E^{\sigma},G^{\sigma})$
structure on $F^{\sigma}$, where $G^{\sigma}=\sigma
G{\sigma}^{-1}$. Remark that we can have $F^{\sigma}=F$ but the
$(E,G)$ of $F$ are not fixed by $\sigma$.

\medskip

{\bf Applications.}

\medskip

Let $E$ be an algebraic extension of a field $k$, and $G$ be a
subgroup of the Galois group of $E$. We want first to answer the
following question, given an  algebraic  extension $F$ of $k$, is
there an $(E,G)$ structure on $F$.  We will find a non commutative
two cocycle with take $H^2(Et_F,L)$ which is the obstruction to
solving this problem, where $L$ is a sheaf on $Et_F$.

\medskip

Recall that, given an $(E,G)$ structure on the field $F$, is
defined by maps $E_i\rightarrow F_i$ where $E_i$ and $F_i$ are
respectively finite products of separable extensions of $E$ and
$F$.

We can remark that the map $E_i\rightarrow F_i$ gives rise to an
$(E,G)$ structure over $Spec(F_i)$, defined as follow:  the
subfamily of etale covers  $F_k$ of $F$ which are etale cover of
$F_i$ is a covering family of $F_i$.

Now suppose that the characteristic of the field is zero, for
every extension (separable) of $F$, there exists always an
extension (separable) $H$ which contains $F$ and $E$. Thus we can
define an $(E,G)$ structure on $H$ which assigns to every etale
morphism $H\rightarrow K$, the map $E\rightarrow H\rightarrow K$.
The category $C(H,E,G)$ of $(E,G)$ structures over $H$ is not
empty.

we can define a sheaf of categories on $Et_F$ defines as follow:

To every etale covering space $H$ of $F$, we assigns the category
$C(E,H,G)$.

This sheaf of categories is a stack. If the group of automorphisms
of each element $H$ of $C(E,F,G)$ is isomorph to  $L(H)$, where
$L$ is a sheaf on $Et_F$, then this stack is a gerbe. The
obstruction of the existence of an $(E,G)$ structure is then given
by a $2-$Cech cocycle.

\medskip

Now, we will give a procedure to build examples of $(E,G)$
structures.

Consider an algebraic extension  $E$ over $k$, with Galois group
$G$. Let $E_S$ be a separable closure of $E$. Suppose that we have
 a representation $h:Gal(E_S:E)\rightarrow G$, and consider an
 element $c$ of $H^1(Gal(E_S:E),G)$ for this representation.
 With this cocycle, we can built a  flat $L$ $G-$bundle on $Et_E$.
 Now consider a trivialization of this bundle. It is defined by a
  family $(E_i)_{i\in I}$ of separable extensions of $E$.
  We define   the $(E,G)$ structure by
  etale morphism $E_i\rightarrow E_i$ . The coordinates change will be given by the coordinates
  change of the bundle $L$. Obviously the holonomy of this $(E,G)$
  structure is $h$.

  \medskip

  {\bf Proposition.}

  {\it If $c=0$, then the $(E,G)$ structure that we have just
  defined is trivial.}

  \medskip

\bigskip

\centerline{\bf Bibliography.}

1. Deligne, P. Le groupe fondamental de la droite projective moins
trois points. Math. Sci. Res. Inst. Publ, 16.

2. Giraud, J. Cohomologie non ab\'elienne. Springer 1971.

3. Goldman, W. Geometric structures on affine manifolds and
varieties of representations. Cont. Math, 74.

4. S.G.A.1 S\'eminaire dirig\'e par A.Grothendieck, Rev\^etement
\'etales et groupe fondamental. Springer.

5. S.G.A.4 S\'eminaire dirig\'e par M, Artin, A. Grothendieck et
J-L Verdier. Th\'eorie des topos et cohomologie \'etales des
sch\'emas. Spinger 1972

6. Micali, A. Corps de fonctions algebriques, notas de
matematicas. Instituto Matematica Nacional de Pesquisas. 1963

\end{document}